\documentclass[12pt]{amsart}

\usepackage{latexsym}
\usepackage[dvips]{graphics,epsfig}
\newtheorem{thm}{Theorem}[section]

\newtheorem{lem}[thm]{Lemma}
\newtheorem{cor}[thm]{Corollary}

\newtheorem{prp}[thm]{Proposition}

\newtheorem{question}[thm]{Question}

\newcommand{\beq}{\begin{equation}}
\newcommand{\eeq}{\end{equation}}
\newcommand{\beqna}{\begin{eqnarray*}}
\newcommand{\eeqna}{\end{eqnarray*}}
\newcommand{\beqn}{\begin{equation*}}
\newcommand{\eeqn}{\end{equation*}}
\newcommand{\bp}{\begin{proof}}
\newcommand{\ep}{\end{proof}}
\newcommand{\bprop}{\begin{proposition}}
\newcommand{\eprop}{\end{proposition}}
\newcommand{\bt}{\begin{theorem}}
\newcommand{\et}{\end{theorem}}
\newcommand{\bex}{\begin{Example}}
\newcommand{\eex}{\end{Example}}
\newcommand{\bc}{\begin{corollary}}
\newcommand{\ec}{\end{corollary}}
\newcommand{\bl}{\begin{lemma}}
\newcommand{\el}{\end{lemma}}

\begin{document}

\title
[On  counterexamples in questions of unique determination]
{On  counterexamples in questions of unique determination of convex bodies}

\author{Dmitry Ryabogin and Vlad Yaskin}
\address{Department of Mathematics, Kent State University,
Kent, OH 44242, USA} \email{ryabogin@math.kent.edu}
\address{Department of Mathematical and
Statistical Sciences, University of Alberta, Edmonton, Alberta T6G 2G1, Canada} \email{vladyaskin@math.ualberta.ca}
\thanks{First author supported in
part by U.S.~National Science Foundation Grants DMS-0652684 and DMS-1101636. Second author supported in part by NSERC}


\begin{abstract}
We discuss a construction that gives counterexamples to various questions of unique determination of convex bodies.
\end{abstract}

\maketitle

\section{Introduction}

This note is motivated by the following question from the book ``Geometric tomography" by R.~J.~Gardner  \cite[Prob. 7.6, p. 289]{G}.

\begin{question}\label{perim}
 {Let $K$ and $L$ be origin-symmetric convex bodies in $\mathbb R^3$ such that the sections $K\cap H$ and $L\cap H$ have equal perimeters for every 2-dimensional subspace $H$ of $\mathbb R^3$.  Is it true that $K=L$?}
\end{question}

The problem has attracted a lot of attention,   but at the moment it is still open; see \cite{HNRZ}, \cite{R}, \cite{Y} for some partial results.
The question can also be generalized to higher dimensions and  all intermediate intrinsic volumes $V_i$, as well as sections of other  dimensions.

In this note we show that the question has a negative answer in the class of general (not necessarily symmetric) convex bodies.

One of the natural steps in generalizing questions involving sections is to replace sections by slabs (or ``thick sections"). We discuss such problems at the end of the paper.

Throughout the paper we use the following notation.  We write $S^{n-1}$  and $B_2^n$  for the unit sphere and closed unit ball in Euclidean space $\mathbb R^n$, and $G(n,k)$ for the Grassmanian of all $k$-dimensional subspaces of $\mathbb R^n$. The hyperplane through the origin with normal vector  $\xi\in S^{n-1}$ is denoted by $\xi^\perp$. A convex body in $\mathbb R^n$ is a compact convex set with nonempty interior. In this paper we always assume that convex bodies contain the origin in their interiors. The orthogonal projection of a convex body $K$ onto a subspace $V$ is denoted by $K|V$.
If $Q$ is a convex body in an $i$-dimensional subspace $V \subseteq \mathbb R^n$, $1\le i\le n$, then we write $\mathrm{vol}_i(Q)$ for its volume.
For a set $C$ in $\mathbb R^n$, its reflection in the origin is denoted by $-C$. We say that a set $C$ is origin symmetric if $C=-C$.
Two sets $C$ and $D$ in $\mathbb R^n$ are said to be congruent if there is an isometry $T$ of $\mathbb R^n$ such that $T(C)=D$.

For a convex body $K$ in $\mathbb R^n$ its support function $h_K$ is defined by $$h_K(x) = \max \{ \langle x,y\rangle: {y\in K}\}, \quad x\in \mathbb R^n.$$  Clearly, $h_K$ is positively homogeneous of degree 1, and therefore is determined by its values on the unit sphere. It is also easy to see that for any subspace $V$ of $\mathbb R^n$ the support function of $K|V$, as a convex body in $V$, is just the restriction of $h_K$ to $V$.

If $K$ is a compact set containing the origin in its interior and star-shaped with respect to the origin, then its radial function $\rho_K$  is defined by $$\rho_K(x)=\max\{\lambda>0: \lambda x\in K\}, \qquad x\in S^{n-1}.$$

We say that a convex body $K$ is of class $C^\infty_+$ if its boundary is a $C^\infty$ manifold in $\mathbb R^n$ with positive Gaussian curvature at each point.

Let $K$ be a compact convex set in $\mathbb R^n$. Its intrinsic volumes $V_i(K)$, $1\le i\le n$, can be defined via Steiner's formula
$$\mathrm{vol}_n(K+\epsilon B_2^n) = \sum_{i=0}^n \kappa_{n-i} V_i(K) \epsilon^{n-i},$$
where the addition is the Minkowski addition,  $\kappa_{n-i}$ is the volume of the $(n-i)$-dimensional Euclidean ball, and $\epsilon\ge 0$.

In particular, if $K$ is a convex body in $\mathbb R^n$, then $V_n(K)$ is its volume, and $V_{n-1}(K)$ is half the surface area.   For these and other facts about intrinsic volumes we refer the reader to the book \cite{S}.

\section{Results}

The idea of the construction that we use in this note goes back to Gardner and Vol\v{c}i\v{c}   \cite{GV} and Goodey, Schneider and Weil
\cite{GSW}.

First we prove an auxiliary lemma.

\begin{lem}\label{lem}
Let $K$ and $L$ be convex bodies in $\mathbb R^n$ containing the origin in their interiors. The following properties are equivalent.
\begin{enumerate}\label{Lem:equivalence}
\item[i)] $K\setminus L = - (L\setminus K)$,

\item[ii)] $\{h_K(\xi), h_K(-\xi)\}=\{h_L(\xi), h_L(-\xi)\}$, for all $\xi \in S^{n-1}$,

\item[iii)] $\{\rho_K(\xi), \rho_K(-\xi)\}=\{\rho_L(\xi), \rho_L(-\xi)\}$, for all $\xi \in S^{n-1}$.
\end{enumerate}
\end{lem}

\bp
The equivalence of (i) and (ii) is shown in \cite{GSW}.

Let us show that (iii) implies (i). Let $x\in K\setminus L$. Since $K$ and $L$ star-shaped with respect to the origin,  $\rho_K(x/|x|)> \rho_L(x/|x|)$. Now (iii) yields $\rho_K(x/|x|)=\rho_L(-x/|x|)$ and $\rho_K(-x/|x|)=\rho_L(x/|x|)$. Therefore, $\rho_K(-x/|x|)< \rho_L(-x/|x|)$, i.e. $x\in - (L\setminus K)$. We have shown that $K\setminus L\subset - (L\setminus K)$. The other inclusion is obtained in a similar fashion. Thus, $K\setminus L = - (L\setminus K)$.

It remains to show that (i) implies (iii). Let $\xi\in  S^{n-1}$. Consider two cases:   $\rho_K(\xi) =\rho_L(\xi)$ and $\rho_K(\xi) \ne \rho_L(\xi)$. In the first case we claim that    $\rho_K(-\xi) =\rho_L(-\xi)$. Indeed, if the latter is not true, then there exists $x$ on the boundary of $K$ (say), such that $x/|x| = -\xi$ and $\rho_K(-\xi) > \rho_L(-\xi)$. But then $-x\in L\setminus K$, meaning that $\rho_K(\xi)<\rho_L(\xi)$. We reach a contradiction.

In the second case, let $\rho_K(\xi) > \rho_L(\xi)$, say. Consider the ray $\{x\in\mathbb R^n: x = t\xi, \  t\ge 0\}.$ It intersects the boundaries of $L$ and $K$ at some points, which we will denote by $p$ and $q$ correspondingly. Since $K\setminus L = - (L\setminus K)$, the ray $\{x\in\mathbb R^n: x = - t\xi, \  t\ge 0\}$ intersects $L\setminus K$ in the interval $(-p, -q]$. Hence, the points $-p$ and $-q$ are on the boundaries of $L$ and $K$ correspondingly. This implies that $\rho_K(\xi) = \rho_L(-\xi)$ and $\rho_K(-\xi) = \rho_L(\xi)$.

\ep

For the reader's convenience, we provide a short proof of the following proposition. Examples of this type appeared in \cite[Thm 3.3.17, Thm 3.3.18]{G} (see also \cite{GV}).

\begin{prp}\label{bodies}
There exist noncongruent convex bodies $K$ and $L$ satisfying the properties of Lemma \ref{Lem:equivalence}. Moreover,  the bodies can chosen to be $C^\infty_+$ bodies of revolution, or polytopes.
\end{prp}

\bp

1) $K$ and $L$ are $C^\infty_+$ bodies of revolution.

We start with two nonzero infinitely
 smooth  functions   $\phi$ and $\psi$,  supported in the intervals $(1/3-\delta,1/3+\delta)$ and $(2/3-\delta,2/3+\delta)$ correspondingly, where $0<\delta < 1/6$.  Next, for
  $t\in [-1,1]$ and $\epsilon>0$,
 we
define
$$
f(t)=\sqrt{1-t^2}+ \epsilon \phi(t)+ \epsilon \psi(t),
$$
$$
 g(t)=\sqrt{1-t^2}+ \epsilon \phi(t)+ \epsilon \psi(-t).
$$
Finally,  we choose $\epsilon$ sufficiently small to guarantee that $f$ and $g$ are concave on $[-1,1]$.

Now we consider  two convex bodies in $\mathbb R^n$,
$$K =\{ x\in \mathbb R^n: x_1^2+\cdots +x_{n-1}^2 \le f^2(x_n)\}\ \mbox{ and } \ L =\{ x\in \mathbb R^n: x_1^2+\cdots +x_{n-1}^2 \le g^2(x_n)\}.$$
The boundaries of $K$ and $L$ are generated by rotating the graphs of the functions $f$ and $g$ defined above.
Since these generating curves are not congruent, it follows that $K$ and $L$ are not congruent.

By construction, $K$ and $L$ are of class $C^\infty$, and also $K\setminus L = - (L\setminus K)$. Finally, choosing $\epsilon>0$ sufficiently small,  $\|\rho_K-\rho_{B_2^n}\|_{C^2(S^{n-1})}$ and $\|\rho_L-\rho_{B_2^n}\|_{C^2(S^{n-1})}$  can be made as small as we wish. By formula (0.41) from \cite{G} this is enough to ensure that $K$ and $L$ have strictly positive curvature.

2) $K$ and $L$ are polytopes.

Let $a_1$,...,$a_n$ be distinct positive numbers and consider the rectangular box $R=[-a_1,a_1] \times \cdots \times [-a_n,a_n]$.
Take any two vertices $u$ and $v$ of $R$ that are connected by an edge. The idea is that $K$ will be constructed by cutting off vertices $u$ and $v$ from $R$, and $L$ by cutting off vertices $u$ and $-v$.

More precisely,
let $F$ and $H$ be supporting hyperplanes to $R$ at the vertices $u$ and $v$ correspondingly, such that $R\cap F= \{u\}$ and $R\cap H= \{v\}$. Suppose $F$ is given by $\langle x,\xi\rangle = t$ for some $\xi\in S^{n-1}$ and $t>0$, and $H$ is given by $\langle x,\eta\rangle = s$ for some $\eta\in S^{n-1}$ and $s>0$. If $\lambda>0$ is sufficiently small, then the sets $R\cap \{x: \langle x,\xi\rangle \ge t -\lambda\}$ and $R\cap \{x: \langle x,\eta\rangle \ge s -\lambda\}$ do not overlap and contain no vertices of $R$ other than $u$ and $v$ correspondingly. Define
$$K = R \cap   \{x: \langle x,\xi\rangle \le t -\lambda\} \cap   \{x: \langle x,\eta\rangle \le s -\lambda\},$$
$$L = R \cap   \{x: \langle x,\xi\rangle \le t -\lambda\} \cap   \{x: \langle x,\eta\rangle \ge -s+\lambda\}.$$

The bodies $K$ and $L$ are not congruent, since the vertices $u$ and $v$ that we cut off from $R$ belong to the same edge, while this is not true for $u$ and $-v$.

Finally, observe that $K$ and $L$ satisfy the property (i) of Lemma \ref{Lem:equivalence}.

\ep

Lemma \ref{lem} and Proposition \ref{bodies} lead to the following result; see Theorems 3.3.17, 3.3.18, 6.2.18, 6.2.19 in \cite{G}.

 \begin{cor}
There are noncongruent convex bodies  $K$ and $L$   in $\mathbb R^n$ such that
  $\forall \xi \in S^{n-1}$,
\begin{equation*}
\mathrm {vol}_{n-1}(K | \xi^{\perp})=\mathrm{vol}_{n-1}(L |\xi^{\perp})
\end{equation*}
and
\begin{equation*}
\mathrm{vol}_{n-1}(K\cap \xi^{\perp})=\mathrm{vol}_{n-1}(L\cap\xi^{\perp}).
\end{equation*}
\end{cor}

This result remains true if  the $(n-1)$-dimensional volume of the projections is replaced by their $i$-th intrinsic volumes. See Gardner \cite[Thm 3.3.17, Thm 3.3.18]{G} for $i=1, n-1$, and Goodey, Schneider and Weil
\cite{GSW} for all other cases.

Now we return to Question \ref{perim} stated in the introduction. The next theorem shows that without the assumption of origin-symmetry, this question would have a negative answer.

\begin{thm}\label{thm1}
There are noncongruent convex bodies $K$ and $L$ in $\mathbb R^n$ such that for all $i$ and $k$ with $1\le i\le k\le n-1$ we have  $$V_i(K\cap H)=V_i(L\cap H) $$
for all  $H\in G(n,k)$.

Moreover, $K$ and $L$ can be constructed in such a way that both of them are  $C^\infty_+$ bodies of revolution, or both $K$ and $L$ are polytopes.

\end{thm}
\bp

Let $K$ and $L$ be noncongruent convex bodies in $\mathbb R^n$ such that $\{\rho_K(\xi), \rho_K(-\xi)\}=\{\rho_L(\xi), \rho_L(-\xi)\}$, for all $\xi \in S^{n-1}$. Recall that   $C^\infty_+$  and polytopal examples of such bodies were given in Proposition \ref{bodies}.

Then, for every $H\in G(n,k)$,  the convex bodies $K\cap H$ and $L\cap H$ in the subspace $H$ satisfy the relation $$\{\rho_{K\cap H}(\xi), \rho_{K\cap H}(-\xi)\}=\{\rho_{L\cap H}(\xi), \rho_{L\cap H}(-\xi)\}, \qquad \forall \xi \in S^{n-1}\cap H.$$ By Lemma \ref{Lem:equivalence}, for these $k$-dimensional bodies we also have $$\{h_{K\cap H}(\xi), h_{K\cap H}(-\xi)\}=\{h_{L\cap H}(\xi), h_{L\cap H}(-\xi)\}, \qquad \forall \xi \in S^{n-1}\cap H,$$  and, in particular, the last relation is true for the projections of $K\cap H$ and $L\cap H$ onto $i$-dimensional subspaces of $H$, $1\le i\le k$. We apply Lemma \ref{Lem:equivalence} to get $\mathrm{vol}_i((K\cap H)|V)= \mathrm{vol}_i((L\cap H)|V)$ for every $i$-dimensional subspace  $V$ of $ H$.

By Kubota's integral recursion (see \cite[p.295]{S}), we have $$V_i(K\cap H) = c_{i,k} \int_{G(H,i)} \mathrm{vol}_i((K\cap H)|V) dV,$$ for some constant $c_{i,k}$. Here, by  $G(H,i)$ we mean the Grassmanian of all $i$-dimensional subspaces of $H$.
Thus, $V_i(K\cap H) =  V_i(K\cap H)$.
 \ep

\smallskip

The following result is proved in \cite{Y}. Let $k$ be an integer with $2\le k\le n-1$. If $K$ and $L$ are origin-symmetric convex polytopes in $\mathbb R^n$ such that $V_{k-1}(K\cap H) = V_{k-1}(K\cap H)$ for all $H\in G(n,k)$, then $K=L$. Theorem \ref{thm1} shows that this result is false without the symmetry assumption.

\smallskip

 It is interesting to note that the following result is true.

\begin{prp}\label{two_sections}
 Let   $i\ne j$ be integers, $1\le i,j\le n-1$, and let $K$ and $L$ be convex bodies (containing the origin in their interiors) such that
 \begin{equation}\label{i-sec}
 \mathrm{vol}_i (K\cap V) = \mathrm{vol}_i (L\cap V), \qquad \forall V\in G(n,i),
 \end{equation}
and
\begin{equation}\label{j-sec}
\mathrm{vol}_j (K\cap W) = \mathrm{vol}_j (L\cap W), \qquad \forall W\in G(n,j).
\end{equation}
Then for all $k$, $1\le k\le n-1$, we have
\begin{equation}\label{k-proj}
\mathrm{vol}_k (K|H) = \mathrm{vol}_k (L|H), \qquad \forall H\in G(n,k).
\end{equation}
\end{prp}
\smallskip

\bp
Conditions (\ref{i-sec}) and (\ref{j-sec}) yield that $\{\rho_K(\xi), \rho_K(-\xi)\}=\{\rho_L(\xi), \rho_L(-\xi)\}$, for all $\xi \in S^{n-1}$; see \cite[Thm 7.2.3, Thm 6.2.16]{G}. By Lemma \ref{Lem:equivalence}, this is equivalent to the condition $\{h_K(\xi), h_K(-\xi)\}=\{h_L(\xi), h_L(-\xi)\}$, for all $\xi \in S^{n-1}$. As shown in  \cite{GSW}, this implies (\ref{k-proj}).

\ep

Clearly, if $K$ and $L$ are two convex bodies such that all their corresponding $i$-dimensional projections have the same $i$-volume, and all corresponding $j$-dimensional projections have the same $j$-volume, then one can not expect that the sections of $K$ and $L$ will have equal volumes (simply because sections are not invariant under translations). Under the same conditions, is it possible to translate $K$ and $L$ so that all their sections would have the same volumes? The answer is again negative. Goodey, Schneider and Weil \cite[p. 86]{GSW} remark that there are convex polytopes $K$  and $L$ such that (\ref{k-proj}) holds for all $k$, but none of their translates satisfy  $K\setminus L = - (L\setminus K)$. Thus, as the proof of Proposition \ref{two_sections} shows, no translates of $K$ and $L$ can satisfy (\ref{i-sec}) and (\ref{j-sec}) for some $i\ne j$.

\smallskip

Our next question concerns slabs in convex bodies.
For $t>0$ and $\xi\in S^{n-1}$ consider the slab of width $2t$ perpendicular to  $\xi$ defined by $$S_t(\xi) = \{x\in \mathbb R^n: |\langle x,\xi \rangle| \le t\}.$$ There are various problems related to volumes of slabs in convex bodies; see e.g. \cite{BK}, \cite{KK}. Here we suggest the following.

\smallskip

\begin{question} Let $K$ and $L$ be origin-symmetric convex bodies in $\mathbb R^n$ that contain the Euclidean ball of radius $t$ in their interiors. Suppose that for some $i$ ($1\le i\le n$) $$V_i(K\cap S_t(\xi) ) = V_i(L\cap S_t(\xi) ), \qquad \forall \xi \in S^{n-1}.$$ Is it true that $K = L$?
\end{question}

Our next result shows that the previous question has a negative answer if we drop the assumption of origin-symmetry.

\begin{thm}\label{thm2}
There are noncongruent convex bodies $K$ and $L$  in $\mathbb R^n$ that contain the Euclidean ball of radius $t>0$ in their interiors  and such that $$V_i(K\cap S_t(\xi) ) = V_i(L\cap S_t(\xi) ), \qquad \forall \xi \in S^{n-1},$$
for all $1\le i\le n$.

Moreover, both $K$ and $L$ can be chosen to be $C^\infty_+$ bodies of revolution or polytopes, as in Theorem \ref{thm1}.

\end{thm}

\bp
Consider the same pair of  bodies $K$ and $L$ as used in Theorem \ref{thm1}. Without loss of generality, we may assume that $t>0$ is small enough so that $tB_2^n$ is contained in the interior of both $K$ and $L$. Since $S_t(\xi)= -S_t(\xi)$, it follows that  $$\{\rho_{K\cap S}(\theta), \rho_{K\cap S}(-\theta)\}=\{\rho_{L\cap S}(\theta), \rho_{L\cap S}(-\theta)\}, \qquad \forall \theta\in S^{n-1},$$ where for brevity we write $S = S_t(\xi)$.
Now, as in Theorem \ref{thm1}, we see that all $i$-dimensional projections of $K\cap S$ and $L\cap S$ have the same  $i$-volumes, and Kubota's integral recursion finishes the proof.
\ep

Since intrinsic volumes are continuous (in the Hausdorff metric) valuations on convex sets, it follows that Theorem \ref{thm1} can be derived from Theorem \ref{thm2}. It is enough to notice that $K\cap S_t(\xi)$ approaches $K \cap \xi^\perp$ in the Hausdorff metric.

\end{document}